\documentclass[a4paper,11pt]{amsart}
\usepackage[T1]{fontenc}
\usepackage[english]{babel}
\usepackage{amscd,graphicx,amsmath,stmaryrd}
\usepackage[arrow,curve,matrix]{xy}
\usepackage{longtable}
 \usepackage{amssymb}

 \def\PV{{ \mathbb{P} }}
 \def\P1{{ \mathbb{P}^1 }}
\def\Pr{{ \mathbb{P}^r }}
\def\Pa{{ {\mathbb{P}^r}^*_a }}

 \def\GG{{\mathcal{G} }}

\def\Z{{\mathbb{Z}}}
\def\Q{{\mathbb{Q}}}

\def\C{{\mathbb{C}}}

\def\hypn{{  \mathcal{H}_{{\Pr}}^{ n}    }}
\def\hypm{{  \mathcal{H}_{{\Pr}}^{ m}    }}

\def\1an{{  \{1 \ldots n \}    }}

\def\unar{{  \{1 \ldots r \}    }}
\def\unam{{  \{1 \ldots m \}    }}
\def\fo{{\mathcal{F}}}

\def\dj{{ \Delta^{J_1}\ldots \Delta^{J_r} O^{ \{ b_1 \} }_{o_1} \ldots O^{ \{ b_s \} }_{o_s}      }}
\def\di{{ \Delta^{I_1}\ldots \Delta^{I_r} O^{A_1}_{o_1} \ldots O^{A_s}_{o_s}      }}
\def\psip{{  \Psi_{\mathcal{P}}   }}

\def\cqfd{\ifmmode\sqw\else{\ifhmode\unskip\fi\nobreak\hfil
\penalty50\hskip1em\null\nobreak\hfil\sqw
\parfillskip=0pt\finalhyphendemerits=0\endgraf}\fi}
\def\sqw{\hbox{\rlap{\leavevmode\raise.3ex\hbox{$\sqcap$}}$%
\sqcup$}}

\newtheorem{theo}{Theorem}
\newtheorem{lem}{Lemma}
\title{}
\author{}

\begin{document}

\title{Algebraic cycles on the Jacobian of a curve with a $g^r_d$}

\author{Fabien Herbaut}

\address{Laboratoire J. A. Dieudonn\'e, Universit\'e
de Nice Sophia-Antipolis, Parc Valrose, 06108 Nice Cedex 2 (France).}

\email{herbaut\char`\@math.unice.fr}

\begin{abstract}
We present relations between cycles with rational coefficients modulo algebraic equivalence on the Jacobian of a curve. These relations depend on the linear systems the curve admits. They are obtained in the tautological ring, the smallest subspace containing (an embedding of) the curve and closed under the basic operations of intersection, Pontryagin product and the pullback and pushdown induced by homotheties.
\end{abstract}

\maketitle

\section{Introduction}\label{int}
In this paper we extend a theorem of Elisabeta Colombo and Bert van Geemen for $d$-gonal curves to linear systems of higher dimension. To express our results in 1.3 we have to recall the Beauville decomposition theorem 1.1 which enlightens the structure of the tautological ring described in 1.2. The main points of the proof are summarised in 1.4.\\ \\
\subsection{} We will use the notation $CH(X)$ for the Chow ring associated to a smooth algebraic variety $X$, and $A(X)$ for its quotient modulo algebraic equivalence. These rings will always be considered \mbox{\textit{tensored by} $\mathit{\Q}$}. The brackets $[V]$ mean we consider the class of a subvariety $V$ of $X$. Let $X$ be an Abelian variety over $\C$ whose group law we denote $m$ and dual variety we denote $\hat{X}$. It admits homotheties $ x \mapsto kx$ which we will also denote $k$ for $k \in \Z$.  We have another product between algebraic cycles, the Pontryagin product, defined between two cycles  $\alpha$ and $\beta'$ by $\alpha*\beta=m_*(p^*\alpha.q^*\beta)$ with $p$ and $q$ the two projections associated to the product $X \times X$. Beauville introduced in \cite{arnauda} the Fourier transform $\fo$,   a $\Q$-linear automorphism between $A(X)$ and $A(\hat{X})$ which exchanges the two products.  He used it to prove in \cite{arnaudb} the decomposition theorem which states that the operators $k_*$ and $k^*$ simultaneously diagonalize.  For each codimension $p$ the subspace $A^p(X)$ splits: \\
\begin{equation}\label{decomp_beauville}  \begin{array}{cccl}
A^p(X)=\bigoplus_{i=p-g}^{g} \ A^p_{(i)}(X)   & \textrm{ where } \alpha \in A^p_{(i)}(X) &\textrm{iff for all } k \in \Z &   k_* \alpha = k^{2g-2p+i} \alpha \\
                            &  &  \big{(} \textrm{or equivalently}               & k^* \alpha = k^{2p-i} \alpha \big{)} \\
& & &
\end{array}
\end{equation}
It is conjectured that $A^p_{(i)}(X)$ vanishes for $i<0$; it is proven if $p \in \{0,1,g-2,g-1,g \}$ .
\subsection{}  Let $C$ be a smooth projective complex curve of genus $g \geq 2$ and $JC$ its Jacobian variety.  As $C$ embeds into $JC$, we can consider the associated cycle $[C]$ in $A^{g-1}(JC)$ and its decomposition\footnote{Perhaps it is useful to recall that $C_{(i)}$ is denoted $\pi_{2g-2-i}C$ in \cite{cite_cvg} and that Polishchuk worked with the classes $p_i=\fo(C_{(i-1)})$ in \cite{poli}.}:
\begin{equation}\label{decomposition_C}
 [C] \ = \ C_{(0)} + \ldots + C_{(g-1)} \ \ \ \ \ \textrm{ where } C_{(i)} \in A^{g-1}_{(i)}(JC)
\end{equation}
\noindent Notice that the cycle $[C]$ doesn't depend on the embedding of the curve, as we work modulo algebraic equivalence. The tautological ring $R$ is the smallest subspace of $A(JC)$ which contains $[C]$ and is stable under the two products and the operators $k_*$ and $k^*$. Actually, $R$ is finite dimensional: Beauville proved in \cite{arnaudc} that $R$ is the subalgebra for the Pontryagin product generated by $C_{(0)}, \ldots C_{(g-1)}$. So we can describe $R$ as a quotient of the $\Q$-algebra $\Q[C_{(0)},\ldots,C_{(g-1)}]$ by an ideal of relations we have to determine. But to discuss the vanishing of algebraic cycles is always a difficult problem. Few results are known for $R$. The cycle $C_{(0)}$ is never zero. Neither are $C_{(1)}$ and $C_{(2)}$ for  generic curves of genus $g\geq 3$ and $g \geq 11$ as Ceresa and Fakhruddin proved in \cite{ceresa} and \cite{fakh} respectively. More recently, Ikeda determined for each degree $d$ a smooth plane curve for which the component $C_{(i)}$ is nonzero in $CH^{g-1}_{(i)}(JC)$ when $i \geq d-3$.
\par
On the other hand, one could ask for vanishing  results. Polishchuk has presented in \cite{poli} an ideal $I_g$ of relations which hold for all curves of genus $g$. For nongeneric curves, Colombo and van Geemen gave the main result. In \cite{cite_cvg}, they stated that for a $d$-gonal curve, the cycle $C_{(i)}$ is algebraically equivalent to zero for $i \geq d-1$. \\ \\
\subsection{} In this paper, we extend this theorem to curves which admit linear systems of higher dimension.  More precisely, we state the following result where we use the abbreviation $g^r_d$ for linear systems of dimension $r$ and degree $d$:
\begin{theo}\label{relations_intro}
Let $C$ be a curve which admits a base point free $g^r_d$. For each integer $s\geq0$ the following relation holds in $A^{g-r}_{(s)}(JC)$ :
\[ \sum_{\scriptstyle 0 \leq a_1,\ldots ,a_r \atop \scriptstyle a_1+\ldots + a_r=s } \beta(d,a_1+1,\ldots,a_r+1) \ C_{(a_1)}* \ldots * C_{(a_r)} = 0 \]
\[ \textrm { \ \ where \ \ \ \ }\beta(d,a_1,\ldots,a_r)=\sum_{i_1=1}^{d} \ldots \sum_{i_r=1}^{d} (-1)^{i_1+\ldots +i_r } \binom{d}{i_1+\ldots + i_r} {i_1}^{a_1} \ldots {i_r}^{a_r} \textrm{ .}\]
\end{theo}
In section \ref{app} we study applications of this theorem to plane and space curves. For almost all genera $g$ there exists a curve with a $g^2_d$ (respectively a $g^3_d$) giving  new relations. By new relations we mean that  they could not be deduced from the $g^1_{d'}$ that the $g^2_d$ (respectively the $g^3_d$) induces and from  knowledge of $I_g$. We list such relations for $g \leq 9$ in the table \ref{tab_plane_curves} and \ref{tab_space_curves}. \\ \\
\indent Do we obtain new algebraic equivalences $C_{(i)}=0$? No, because the monomial relations we obtain are  $B(r,d,g) \ C_{(d-2r+1)} \ = \ 0$, where for each dimension $r$ the integer $B(r,d,g)$ is the  number of $(r-2)$ planes which cut the curve (mapped to $\Pr$ by the $g^r_d$) in $2r-2$ points. When such a situation arises we can  construct (by projection) a $g^1_{d-2r+2}$. In this case, the monomial relation above does not teach us more than Colombo and van Geemen's theorem. We explain it in section \nolinebreak \ref{sec}. \\ \\
\par
\subsection{} The first sections are devoted to the proof of theorem \ref{relations_intro}. Let $C_d$ be the $d$-fold symmetric product of the curve $C$. Choosing an element of  $C_d$ provides a morphism \mbox{$u_d \ : \ C_d \longrightarrow JC$}. Throughout the whole paper, $\GG$ will design a $g^r_d$ and $G_n$ the truncated linear system of degree $n$: this is the set of divisors of $C_n$ we can complete to a divisor of $\GG$.
\[  G_n = \{ D \in C_n \mid \exists \  E \in C_{d-n} \ , \  D+E \in \mathcal{G}  \} \ \ \ \ \ \textrm{ for  }1 \leq n\leq d \]
\noindent  These truncated linear systems may be considered as subvarieties\footnote{They have already been used by Kouvidakis in \cite{kouvidakis} to construct curves in the symmetric product whose Neron Severi classes are known. Izadi has also used the truncated  linear systems in \cite{izadi} to study how curves deform in Abelian varieties.} of the symmetric products $C_n$. The scheme structure is described in \cite{acgh}, paragraph 3 of chapter VIII. They play an important role in Colombo and van Geemen's proof.  The cornerstone of the present note is the generalization of a formula\footnote{Proposition (3.4) of \cite{cite_cvg} states that if $\GG$ is $g^1_d$, then we have the equality in $A^{g-1}(JC)$:  \\ \[ {u_n}_*[G_n]=\sum_{i=1}^{n} \frac{(-1)^{i-1}}{i} \binom{d}{n-i} \ i_*C \]} they obtained in case $r=1$ which expresses the algebraic classes $[G_n]$ (and especially $[\GG]=[G_d]$) as functions of the diagonals in $C_d$ :
\begin{theo}\label{th_gk_general_intro}
If $C$ admits a base point free $g^r_d$ then the following equivalences hold in $A^{d-r}(C_d)$ for $r \leq n \leq d$ :
\[ [G_n] \ =   \sum_{ \scriptstyle 1 \leq i_1 \leq \ldots \leq i_r   }  \binom{d}{n-\sum i_u} \Big{(}\prod_{u=1}^{r} \frac{(-1)^{i_u-1}}{i_u} \Big{)} \ \  [\delta_{i_1,\ldots,i_r}+(n-\sum_{u=1}^r i_u)o ]
\]
\noindent where  we choose a point $o$ of $C$ and we write $\delta_{i_1,\ldots,i_r}$  for the generalized diagonals in $C_d$ :
\[ \delta_{i_1,\ldots,i_r}\ = \
\{ i_1x_1+\ldots+i_rx_r  \ \mid \ x_i \in C \}
\]
\end{theo}
\noindent We state in section \ref{induction_relation} a relation between $[G_r],\ldots,[G_n]$ for all $r \leq n \leq d$ . It enables us to prove theorem \ref{th_gk_general_intro} by induction in section \ref{rel}. In section \ref{sec}, we notice that ${u_d}_*[G_d]$ is zero in $A(JC)$, because $G_d$ is $\GG$ and $u_d$ contracts this projective space onto a point in the Jacobian. \\ \\

\textit{Acknowledgements: } This paper is part of my thesis written at the University of Nice Sophia-Antipolis,  in the Laboratoire Jean Dieudonn\'e. I would like to thank my Ph.\ D.\ advisor, Prof. Arnaud Beauville, for his guidance and his constant support. I am  grateful to Claire Voisin for much useful advices. It is also a pleasure to thank George Hitching for his help with the English.

\section{ Brief reminder about Fourier transform and tautological ring.}
The results recalled here are proven in \cite{arnauda} and \cite{arnaudb}.
\subsection{} We will identify $JC$ and its dual to consider the Fourier transform $\fo \ : \ A(JC) \rightarrow A(JC)$. It can be  defined by the correspondence associated to the exponential of the class of the Poincar\'e line bundle in $A(JC \times JC)$. As \mbox{$ \fo \  \circ\ \fo= (-1)^g (-1)_*$},  this is an isomorphism. It exchanges the two products, that is:
\[ \fo(x.y) \ = \ (-1)^g \ \fo(x) * \fo(y) \ \ \textrm{ and } \ \ \fo(x*y) \ = \ \fo{x}.\fo{y}\]
\subsection{} The decomposition (\ref{decomp_beauville}) leads to a bigraduation, in the sense that:
\[ A^p_{(i)} . A^q_{(j)} \subset A^{p+q}_{(i+j)} \]
\noindent The proof of the existence of the decomposition (\ref{decomp_beauville}) gives:
\[ \fo{A^{p}_{(i)}} \ = \ A^{g-p+i}_{(i)}\]
\noindent and so the Pontryagin product is also homogeneous of degree $-g$ respect to  the bigraduation :
\[ A^p_{(i)} *  A^q_{(j)} \subset A^{p+q-g}_{(i+j)} \]

\subsection{} The subspace $R$ is spanned by homogeneous elements for both graduations, so $R$ is \mbox{bigraded} and if we note $R^p_ {(i)}$ for the intersection $R^p_ {(i)} =R \cap A^p_{(i)}(JC)$, we have :
\[  R = \bigoplus_{{ \scriptstyle 0 \leq p \leq g \atop \scriptstyle 0 \leq i \leq p } } R^{p}_{(i)}\]

\section{Relations between truncated linear systems}\label{rec}
Let us consider $\GG$ a base point free $g^r_d$. It induces a nondegenerate morphism\footnote{This amounts to saying that no hyperplane of $\Pr$ contains  $\Phi(C)$.} $\Phi \ : \ C \rightarrow \Pr$  and truncated linear systems $G_n$. We will also consider for positive integers $r,n$ the subvariety $\hypn \subset (\Pr)^n$ of $n$-tuples whose components are contained in a hyperplane:
\begin{equation}
\hypn = \{ (y_i)\in (\Pr)^n \ | \ \exists H \in \Pr^* \ \forall i \ y_i \in H   \}
\end{equation}
\noindent Now we introduce the morphisms $\Psi_{\mathcal{P}} \ : \ C^k \rightarrow C^n$ which are defined for integers $0 \leq k \leq n$ and for $\mathcal{P}$ an ordered $k$-partition of $\1an$ (this is a partition of $\1an$ into $k$ ordered subsets $A_1,\ldots,A_k$). We set :
\[ \begin{array}{ccccc}
\Psi_{(A_1,\ldots,A_k)} & : & C^k & \longrightarrow& C^{n} \\
                                                &  & x      & \longmapsto &  y
\end{array} \]
\noindent where $y_i=x_j$ if and only if $j$ is the only integer such that $i \in A_j$. For example, we have $\Psi_{ \{1,2\}}  :  C  \rightarrow C^{2}$ defined by $x \mapsto (x,x)$. Besides, $\Psi_{(\{1,2 \},\{ 3 \})} \ : \ C^2  \rightarrow C^{3}$ is defined by $ (x,y)\mapsto (x,x,y)$. We will also note $\sigma_k$ for the addition morphism : $\sigma_k : C^k \rightarrow C_k$. Now we can state the main result of this section:
\begin{theo}\label{rel_rec_classes_Gk}
We have for each integer $n \geq r$ the equality in $CH(C^n)$
\[ {\Phi^n}^* [\mathcal{H}^{n}_{\Pr}] =  \sum_{k=r}^{n} \frac{1}{k!}\sum_{\scriptstyle  \mathcal{P} \textrm{ ordered }\atop \scriptstyle \textrm{ } k- \textrm{partition} \textrm{ of }  \1an } {\Psi_{\mathcal{P}}}_* ({\sigma_k}^* [G_k])\]
\end{theo} \ \\
We first state a set theoretic analogue of theorem \ref{rel_rec_classes_Gk} in \ref{set_theoric}. We conclude in \ref{transversality} when we prove that $(\Phi^n)^{-1} \hypn$ is reduced. In \ref{hnpr} we compute the class of $\hypn$ in $CH((\Pr)^n)$.

\subsection{}\label{hnpr} One has, by the proposition 8.4 in \cite{fulton}, that
\begin{equation}
CH(\Pr)=\frac{\Q[h]}{(h^{r+1})}
\end{equation}
\noindent where $h$ is the class of a hyperplane of $\Pr$. By proposition 8.3.7 in \cite{fulton} we have
$$ CH((\Pr)^n) \simeq  \frac{\Q[h_1,\ldots,h_n]}{(h_1^{r+1},\ldots,h_n^{r+1})} $$
where $h_i$ is the class of an hyperplane in $\Pr$  at the position $i$. Following this isomorphism we can state
\begin{theo}\label{class_hnpr}
For integers $0 \leq r \leq n$ we have,  in $CH((\Pr)^n)$,
$$[\hypn] = \sum_{  \scriptstyle I \subset \{ 1,\ldots,n \}
                             \atop \scriptstyle \# I=n-r                  } \Big{(}\prod_{a \in I} h_a \Big{)} \textrm{ .}$$
\end{theo}
\noindent I would like to thank Claire Voisin for the following proof.
\begin{proof}Let us consider the incidence variety $\mathcal{I} = \{ (x,H) \in \Pr  \times {\Pr}^* \ \mid \ x \in H\} $. Let $\mathcal{B}$ be a basis of $\C^{r+1}$, and $\mathcal{B}^*$ the dual basis. Such bases induce systems of coordinates on $\Pr$ and $(\Pr)^*$. If $([x_0:\ldots:x_r],[a_0:\ldots a_r])$ represents an element of $\Pr \times {\Pr}^*$, the equation of $\mathcal{I}$ is $\sum a_ix_i=0$. So $\mathcal{I}$ is a hypersurface of bidegree $(1,1)$ in $\Pr \times {\Pr}^*$. Now let us consider the variety $ I= \{ (x_1,\ldots,x_n,H) \ \mid \ \forall i \in \1an \ x_i \in H \} $.

$$ \xymatrix{
 & I \ \ \ \subset  (\Pr)^n  \times  \Pr^*  \ar[ldd]_{\pi} \ar[rdd]^{p_i}&    \\
 &                   &                    \\
(\Pr)^n &                                    & \Pr \times \Pr^* \\
} $$

\noindent The projection $p_i$ maps $(x_1,\ldots,x_n,H)$ onto $(x_i,H)$ while $\pi$ maps it onto $(x_1,\ldots,x_n)$. We have $\hypn=\pi(I)$, and $I$ is the transverse intersection $\bigcap_{i=1}^n \ p_i^{-1}(\mathcal{I})$. Following the isomorphism
$$ CH((\Pr)^n \times \Pr ^*) \simeq  \frac{\Q[h_1,\ldots,h_n,h]}{(h_1^{r+1},\ldots,h_n^{r+1},h^{r+1})} $$
where $h_i$ is the class of a hyperplane at position $i$ and $h$ the class of a hyperplane in $\Pr^*$, we have in $CH^{n}((\Pr)^n \times \Pr ^*)$ the equality
\begin{equation}\label{product}
 [I] \ = \ \prod_{i=1}^n \ \ (h_i+h) \ \ \ \ \textrm{ .}
\end{equation}
Lastly, taking the pushdown $\pi_*$ is the same as considering the coefficient of $h^r$ in the product (\ref{product}). \end{proof}

\subsection{}\label{set_theoric} Here we establish a set theoretic equality.
\begin{theo}\label{rel_rec_set}
For integers $0 \leq r \leq n$ we have
$$ ({\Phi^n})^{-1} \hypn = \bigcup_{k=r}^{n} \bigcup_{ \scriptstyle  \mathcal{P} \textrm{ } k- \textrm{partition} \atop \scriptstyle \textrm{ of }\1an  } \Psi_{\mathcal{P}}(\sigma_k^{-1}G_k) \textrm{.}$$
\end{theo}
\begin{proof}Let us prove the inclusion of  the left hand side in the right one (the other inclusion is straightforward). Let $x$ be an element of $({\Phi^n})^{-1} \hypn $. For integers $i,j \in \1an$ we will denote $i \sim j$ if and only if $x_i=x_j$ for generic $x$. It defines an equivalence relation, so we deduce a partition $\mathcal{P}$ of $\1an$. Let us choose an order : $\mathcal{P}=(A_1,\ldots,A_k)$. If $\alpha_1,\ldots,\alpha_k$ are representatives of these $k$ classes, the components $x_{\alpha_1},\ldots,x_{\alpha_k}$ are generically distinct and mapped  by $\Phi$ into a hyperplane. So $(x_{\alpha_1},\ldots,x_{\alpha_1})$ is an element of $\sigma_k^{-1}G_k$, and $x$ is an element of $\Psi_{\mathcal{P}}(\sigma_k^{-1}G_k)$. We thus have :
$$ {(\Phi^n)}^{-1} \hypn = \bigcup_{k=1}^{n} \bigcup_{ \scriptstyle  \mathcal{P} \textrm{ } k- \textrm{partition} \atop \scriptstyle \textrm{ of }\1an  } \Psi_{\mathcal{P}}({\sigma_k}^{-1}G_k) $$
But for $k<r$, we know that $k$ points of $\Pr$ always lie in a hyperplane, so $G_k=C_k$ and $\sigma_k^{-1}G_k=C^k$. The theorem follows.
\[  \bigcup_{ \scriptstyle  \mathcal{P} \textrm{ } k- \textrm{partition} \atop \scriptstyle \textrm{ of }\1an  }
 \Psi_{\mathcal{P}} (C^k) \subset
 \bigcup_{ \scriptstyle  \mathcal{P} \textrm{ } r- \textrm{partition} \atop \scriptstyle \textrm{ of }\1an  }
  \Psi_{\mathcal{P}} (C^r) \]
\end{proof}
\subsection{}\label{lemmas_prel}
We regroup in the following theorem the preliminary lemmas we need to prove that $(\Phi^n)^{-1} \hypn$ is reduced:
\begin{theo}\label{prelimilemmes}
Let $\Gamma$ be the intersection $\Phi(C)^n \cap \hypn$. We have : \\ \\
a) $\Gamma$ is pure of dimension $r$. \\
b) For $p=(p_1,\ldots,p_n)$ a generic point of $\Gamma$, the space $<p_1,\ldots,p_n>$ spanned by the components $p_i$  is of dimension $r-1$. \\
c) For $p=(p_1,\ldots,p_n)$ a generic point of $\Gamma$, for all $i \in \1an$ the projective tangent space $\mathbb{T}_{p_i}\Phi(C)$ cuts the projective space $<p_1,\ldots,p_n>$ transversally. \\
d) Let $V$ be an irreducible component of $\Gamma$. For all $i \in \1an$, the projection $p_i \ : \ (\Pr)^n \rightarrow \Pr$ is not constant on $V$.
\end{theo}

\begin{proof} As $\Gamma$ is the intersection in $(\Pr)^n$  of two projective varieties of dimensions $nr-n-r$ and $n$, its dimension verify $dim(\Gamma) \geq r$. Let us consider the incidence variety
\[  I= \{ (x_1,\ldots,x_n,H) \ \mid \ \forall i \in \1an \ x_i \in H \} \]
\noindent and the projections
$$\xymatrix{
 & I \ \ \ \subset  \Phi(C)^n  \times  \Pr^*  \ar[ldd]_{\pi_1} \ar[rdd]^{\pi_2}&    \\
 &                   &                    \\
\Phi(C)^n &                                    &\Pr^*  \\
}$$
The curve $\Phi(C)$ is nondegenerate, so the fiber above a hyperplane $H$ is 0-dimensional. It enables us to bound the dimension of each irreducible component of $I$ by $r$, and so, of each irreducible component of $\pi_1(I)=\Gamma$. This proves \textit{a)}. The proposition \textit{b)} is true when $r$ equals $1$. When $r$ is greater than $1$, let us consider an irreducible component $V$  of $I$. We can find $r$ pairwise distinct components on an open set $U \subset V$, otherwise we could bound the dimension of $V$ by $r-1$. Now recall the general position theorem as stated in the chapter 3 of \cite{acgh}. \\ \\
\noindent \textbf{General position Theorem: } Let $\mathcal{C} \subset \Pr$, $r \geq 2$, be an irreducible nondegenerate possibly singular curve of degree $d$. Then a general hyperplane meets $\mathcal{C}$ in $d$ points, any $r$ of which are linearly independent. \\ \\
Let $W \subset \Pr^*$ be the open set whose points correspond to such hyperplanes. Then the $r$ generically distinct components we chose above are linearly independent on the open set $U \cap \pi_2^{-1}(W)$, which proves \textit{b)}. The hyperplanes of $U \cap \pi_2^{-1}(W)$ cut $\phi(C)$ transversally because they cut $\Phi(C)$ in $d$ points, so \textit{c)} is proven. To prove \textit{d)}, let us suppose that a component is constant on $V$, for example the first of them equals $a$. We note $\Pa$ for the hyperplanes of $\Pr$ which contain $a$, and $I_a$ for the incidence variety $I_{a}=\{ (a,y_2\ldots,y_n,H) \in \{ a \} \times \Phi(C)^{n-1} \times \Pr^*_{a} \ \mid \ \forall i \ y_i\in H \} $. It admits projections :
$$\xymatrix{
 & I_{a} \ \ \subset  \{ a \} \times \Phi(C)^{n-1}  \times   \Pa \ar[ldd]_{q_1} \ar[rdd]^{q_2} &    \\
 &           &      \\
\{ a \} \times \Phi(C)^{n-1} &      &\Pa  \\
}$$
The dimension of each irreducible component of $I_a$  is bounded by $r-1$, because the fiber above a point of $\Pa$ is 0-dimensional. As $V$ is the image by $q_1$ of an irreducible component of $I_a$, its dimension is also bounded by $r-1$.
\end{proof}
\subsection{}\label{transversality}
 What is the tangent space to the variety $\hypn$ at a point $p=(p_1,\ldots,p_n)$ of $(\Pr)^n$?  As the $p_i$ lie in a same hyperplane, we can choose a system of coordinates such that the points $p_i$ can be expressed as the columns of the matrix  \\ \\

$$\begin{pmatrix}
1 & & 1 \\
v^1_1 &  & v^n_1\\
\ldots & & \ldots \\
\ldots & \ldots & \ldots \\
\ldots & & \ldots \\
v^1_{r-1} & & v^n_{r-1} \\
0  & & 0\\
\end{pmatrix} $$\ \\ \\

\noindent In a neighbourhood of $(p_1,\ldots,p_n)$, the coordinates of a point of $(\Pr)^n$ are  \\ \\
\begin{equation}
\begin{pmatrix}
1 & & 1 \\
v^1_1+ \epsilon^{1}_{1} &  & v_1^n + \epsilon^{n}_{1}\\
\ldots & & \ldots \\
\ldots & \ldots & \ldots \\
\ldots & & \ldots \\
v^1_{r-1}+ \epsilon^{1}_{r-1}  & & v^n_{r-1} + \epsilon^{n}_{r-1}\\
\epsilon^1_r  & & \epsilon^n_r\\
\end{pmatrix}
\end{equation} \ \\

It corresponds to a point of $\hypn$ if and only if all $(r+1) \times (r+1)$ minors of the above matrix vanish. The terms of degree one in $\epsilon^i_j$ of these minors are the determinants of the following matrix  \\ \\
\begin{equation}\label{points_modifies}
\begin{vmatrix}
1 & \ldots &\ldots &\ldots &\ldots &\ldots &1 \\
v_1^{\sigma(1)} &  \ldots & \ldots & \ldots & \ldots & \ldots & v_1^{\sigma(r+1)} \\
& & & \ldots  & & & \\
& & & \ldots  & & & \\
v_{r-1}^{\sigma(1)} &  \ldots & \ldots & \ldots & \ldots & \ldots & v_{r-1}^{\sigma(r+1)} \\
\epsilon_{r}^{\sigma(1)}  &  \ldots & \ldots & \ldots & \ldots & \ldots & \epsilon_{r}^{\sigma(r+1)}
\end{vmatrix}
\end{equation} \ \\
\noindent for each injective maps $\sigma : \{ 1\ldots r+1 \} \rightarrow \1an$ . So the tangent vectors correspond to the vectors  $(\epsilon_j^i)_{\scriptstyle  i \in \{ 1 \ldots n\}\atop \scriptstyle j \in \{ 1 \ldots r\}}$ of $\C^{nr}$ with conditions on the  $\epsilon^i_{r+1}$: if we choose $r+1$ points $p_{\sigma(1)}, \ldots,p_{\sigma(r+1)}$ among $p_1,\ldots,p_n$, the  points which \mbox{correspond} to the columns  of  \ref{points_modifies} should lie in a hyperplane. \\
\indent Note that if the points $p_1,\ldots,p_n$ span a projective space of dimension strictly less than $r-1$, this last  condition is always satisfied. So the tangent space $T_p \hypn$ is all of $T_p (\Pr)^n$, and  $p$ is singular in $\hypn$. \\
\indent Now  suppose that these points span a projective space of maximal dimension, for example the first $r$  are free. For all $j \in \{ r+1 \ldots n\}$, the following determinant vanishes:\\

$$\begin{vmatrix}
1 & \ldots &\ldots &\ldots &\ldots &1 &1 \\
v_1^{1} &  \ldots & \ldots & \ldots & \ldots & v^r_1 & v_1^j \\
& & & \ldots  & & & \\
& & & \ldots  & & & \\
v_{r-1}^{1} &  \ldots & \ldots & \ldots & \ldots & v^r_{r-1} & v^j_{r-1} \\
\epsilon_{r}^{1}  &  \ldots & \ldots & \ldots & \ldots & \epsilon^{r}_{r} & \epsilon^j_r \\
\end{vmatrix}$$ \ \\ \\
So we have a relation between $\epsilon^1_r,\ldots,\epsilon^r_{r}$ and $\epsilon^j_r$. The coefficient  before $\epsilon^j_r$ is not zero, because it is the top left  $(r-1) \times (r-1)$ minor. We can thus determine $n-r$ free relations, so the dimension of the tangent space is less than $nr-(n-r)$. But this is the dimension of  $\hypn$, so the dimensions are equal and $p$ is nonsingular in $\hypn$. \\ \\
\indent Nonsingular points of $\hypn$ are exactly the points $p=(p_1,\ldots,p_n)$ whose components $p_i$ span a projective space of dimension $r-1$.  According to part b) of theorem \ref{prelimilemmes}, a generic point of $\Gamma$ is nonsingular on $\hypn$. \\ \\
\indent A point $(p_1,\ldots,p_n) \in \Phi(C)^n$ is singular if and only if one component  $p_i$ is singular in $\Phi(C)$. By part c) of the theorem  \ref{prelimilemmes}, a generic point of $\Gamma$ is nonsingular in $\Phi(C)^n$. \\ \\
\indent Now let $p$ be a nonsingular point of $\Phi(C)^n$ and $\hypn$. By the computation of the tangent space at $\hypn$ in $p$, the sum of the tangent spaces in $p$ of the two subvarieties is $T_p(\Pr)^n$ as soon as for every $i \in \1an$ the projective tangent space  $\mathbb{T}_{p_i}\Phi(C)$ cut transversely the space spanned by the components $p_i$. This is the case for generic $p$ in $\Gamma$ by part d) of theorem  \ref{prelimilemmes}. \\ \\
Now, the morphism $\Phi^n : C^n \rightarrow \Phi(C)^n $ is a finite morphism. Its ramification locus corresponds to the  $n$-tuple $(p_1,\ldots,p_n)$ such that one of the components $p_i$ is in the ramification locus of  $\Phi$. According to part a) of theorem \ref{prelimilemmes}, no irreducible component of $\Phi(C)^n \cap \hypn$ is contained in the ramification locus of  $\Phi(C)^n$.  We conclude that the schematic preimage of the intersection is reduced. We can thus deduce theorem \ref{rel_rec_classes_Gk} from theorem \ref{rel_rec_set}.

\section{Classes of the truncated linear systems}\label{induction_relation}
We again consider $\GG$, a base point free $ g^r_d$ . We fix a divisor $D=p_1 + \ldots + p_d$ of $\GG$. We give below the classes of the truncated linear systems $[G_k]$ and $\sigma_k^*[G_k]$ in the Chow rings $CH^{g-r}(C_k)$ and $CH^{g-r}(C^k)$. In the following theorem, the sums are taken   \\ \\
\noindent - (for the first one) over the partitions $I_1$, \ldots, $I_r$ of $\{ 1,\ldots,n\}$  we don't order, except by growing cardinals. So if we note $i_u$ the cardinal of $I_u$, we shall have $i_1 \leq \ldots \leq i_r$. For each choice of  $I_1, \ldots, I_r$, we note $a_1, \ldots, a_{n-\sum i_u}$ the elements of $\1an \setminus \bigcup I_k$ we order from the smallest to the largest.\\ \\

\noindent - the $n- \sum i_u$ distinct points $o_1, \ldots, o_{n-\sum i_u}$ chosen in the support of the divisor $D$. These points are considered ordered in the sum i) and unordered in the sum ii).
\begin{theo}\label{classe_Gk}
If $C$ admits $\GG$ a base point free $g^r_d$ we have the following equalities in  $CH^{g-r}(C^n)$ and  $CH^{g-r}(C_n)$ respectively :
\begin{displaymath}
 \textrm{i)} \ \ \ \sigma_n^*[G_n] =   \sum_{ \scriptstyle I_1,\ldots,I_r \subset \1an \atop \scriptstyle o_1,\ldots,o_{n-\sum i_u} \textrm{ distinct } }  \Big{(}\prod_{u=1}^r (-1)^{i_u-1} (i_u-1)! \Big{)} \ \ [ \ \Delta^{I_1} \ldots \Delta^{I_r} O_{o_1}^{ \{ a_1 \} } \ldots O_{o_{n-\sum i_u}}^{ \{ a_{n-\sum i_u} \} \  } ]
\end{displaymath}
\begin{displaymath}
 \textrm{ii)} \ \ \ [G_n] =   \sum_{ \scriptstyle 1 \leq i_1 \leq \ldots \leq i_r  \atop \scriptstyle o_1,\ldots,o_{n-\sum i_u} \textrm{ distinct }   } \Big{(}\prod_{u=1}^{r} \frac{(-1)^{i_u-1}}{i_u} \Big{)} \ \  [ \ \delta_{i_1,\ldots,i_r}+o_1+\ldots+o_{n-\sum i_u} \ ]
 \end{displaymath}
\end{theo}
\subsection{} We will prove in this subsection that the second proposition is a consequence of the first one.
\noindent  Let us write $s$ for  $n-\sum i_u$. The image by $\sigma_n$ of a generalized diagonal $\Delta^{I_1} \ldots \Delta^{I_r}O_{o_1}^{\{a_1 \} }\ldots O_{o_s}^{ \{a_s \}}$ is $\delta_{i_1,\ldots,i_r}+o_1 + \ldots +o_{s}$. When the cardinals $i_l$ are pairwise distinct  we find above a generic point of  $\delta_{i_1,\ldots,i_r}+o_1 + \ldots +o_{s}$ one preimage  by $\sigma_n$. Else, note $d_1, \ldots d_t$ the integers such that  $i_1 = \ldots = i_{d_1}$, $i_{d_1} \neq i_{d_1+1}$, $i_{d_1 + 1}= \ldots = i_{d_1 + d_2}$, $i_{d_1+d_2} \neq i_{d_1 + d_2 +1}$, \ldots, $i_{d_1+\ldots+d_{t-1}+1}=\ldots=i_{d_1+\ldots+d_t}$. Let $x$ be a generic point of $\delta_{i_1,\ldots,i_r}+o_1 + \ldots +o_{s}$. If $d_1$ elements $x_1, \ldots x_{d_1}$ appear with multiplicity $i_1=\ldots=i_{d_1}$, there are $d_1!$ ways to associate them with the $d_1$ sets  $I_1, \ldots I_{d_1}$. We then count  $d_1! \ldots d_t !$ antecedents of $x$ by $\sigma_n$ in $\Delta^{I_1} \ldots \Delta^{I_r}O_{o_1}^{\{a_1 \} }\ldots O_{o_s}^{ \{a_s \}}$, so
$$ {\sigma_n}_* [\Delta^{I_1} \ldots \Delta^{I_r}O_{o_1}^{\{a_1 \} }\ldots O_{o_s}^{ \{a_s \}}]=d_1! \ldots d_t! \ [\delta_{i_1,\ldots,i_r}+o_1+\ldots+o_s]\textrm{.}$$
Let us count the diagonals $\Delta^{I_1} \ldots \Delta^{I_r}O_{o_1}^{\{a_1 \} }\ldots O_{o_s}^{ \{a_s \}}$ that $\sigma_n$ map to $\delta_{i_1,\ldots,i_r}+o_1+\ldots+o_s$.  Let us count the ways to partition $\1an$ in $d_1$ sets with $i_1$ elements, \ldots, $d_r$ sets with $i_r$ elements and one pointed set with $(n-(i_1 + \ldots +i_r))$ elements. We find
\[
 \frac{1}{d_1! \ldots d_r!} \ \binom{d}{i_1,i_2,\ldots,i_r,n-s} \textrm{ such ways, that is }
 \]
\[
 \frac{1}{d_1! \ldots d_t!} \ \frac{d!}{i_1! \ldots i_r! (n-s)!}
\]
such ways. At last, there exist $(n-s)!$  ways to permute the $n-s$ points of $D$ we have chosen, because for each permutation $\tau$ of $\mathcal{S}_{n-s}$ we have
\[
 \sigma_n \big{(}\Delta^{I_1} \ldots \Delta^{I_r}O_{o_{\tau.1}}^{\{a_1 \} }\ldots O_{o_{\tau.s}}^{ \{a_s \}} \big{)} \ = \ \delta_{i_1,\ldots,i_r}+o_1+\ldots+o_s\textrm{.}
\]
The morphism $\sigma_n$ is degree $n!$ and we conclude with the pullback-pushdown formula:
\[ {\sigma_n}_* \circ {\sigma_n}^{*} =deg(\sigma_n).Id_{CH(C_n)} \textrm{.}\]
\subsection{}
\noindent  We know that $G_r=C_r$ and $\sigma_r^{-1}G_r=C^r$. We deduce that $ \sigma_r^*[G_r]=[C^r]$. This is the theorem \ref{classe_Gk} for $n=r$.
Suppose that theorem \ref{classe_Gk} is proven for $n \leq m-1$. By theorem \ref{rel_rec_classes_Gk}  we have
\begin{equation}\label{audessus}
  [G_m]={\Phi^m}^* [\hypm] - \sum_{k=1}^{m-1} \frac{1}{k!} \sum_{\scriptstyle  \mathcal{P} \textrm{ } \textrm{ordered } k-\textrm{partition} \atop \scriptstyle \textrm{ of }\unam } {\Psi_{\mathcal{P}}}_* ({\sigma_k}^* [G_k])
\end{equation}
and theorem \ref{class_hnpr} gives the class of $\hypm$ in  $CH((\Pr)^m)$
\[
 [\hypm] = \sum_{  \scriptstyle I \subset \{ 1,\ldots,m \}
                             \atop \scriptstyle \# I=m-r                  } \Big{(}\prod_{a \in I} h_a \Big{)}
\]
The divisor $D=p_1+ \ldots + p_d$ is defined as the pullback of the classe of an hyperplane of $\Pr$. So we can write
\begin{equation}\label{margot}
{\Phi^m}^* [\hypm]= \sum_{ \scriptstyle o_1,\ldots,o_{m-r} \in \{ p_1,\ldots,p_d\}
                             \atop \scriptstyle 1 \leq a_1 < \ldots < a_{m-r} \leq m     } O_{o_1}^{ \{a_1\}} \ldots O_{o_r}^{ \{ a_r \} }\textrm{.}
\end{equation}
where the points $o_i$ are chosen in the support of $D$ and are not necessarily distinct. By the induction hypothesis the classes $\sigma_k^*[G_k]$ are sums of classes of varieties $\dj$ where $(J_1,\ldots,J_r,\{ o_1 \},\ldots,\{ o_r \} )$ is a partition of $\{ 1,\ldots,k\}$. By definition of the morphisms $\psip$, the classes
\begin{displaymath}
\psip_*[\dj]
\end{displaymath}
 are again classes of the form:
\[ [\di] \]
\noindent  with $(I_1,\ldots,I_r,A_1,\ldots,A_s )$  a partition of $\{ 1,\ldots,m\}$. We can express the cycle $\sigma_m^*[G_m]$ as a linear combination of such classes. \\ \\

We will fix a  partition $(I_1,\ldots,I_r,A_1,\ldots,A_s )$ of $\{ 1,\ldots,m\}$ and look for which integer $k$ of $\{ r, \ldots, n \}$ and which partition $(J_1,\ldots,J_r,\{ b_1 \},\ldots,\{ b_s \} )$ of $\{1,\ldots,k\}$ there exists a $k$-ordered partition $\mathcal{P}$ of $\unam$ such that
\[ \psip_*  [\dj]=\di \]
\noindent We also have  to determine the coefficients associated to these classes in the sum (\ref{audessus}). \\ \\
Since for each permutation $\sigma \in \mathcal{S}_r$ of the sets  $J_1,\ldots,J_r$ the varieties considered are the same, we deduce that
\[
\dj = { \Delta^{J_{\sigma(1)}}\ldots \Delta^{J_{\sigma(r)}} O^{A_1}_{o_1} \ldots O^{A_s}_{o_s}      }
\]
\noindent so we will consider the ones such that for all $u \in \unam$,
\[ \psip \Delta^{J_u}= \Delta^{I_u}\]
If we denote $i_u$ and  $j_u$ the cardinals of the sets $I_u$ and $J_u$, we  have by definition of the varieties $\Delta$ and the morphisms $\psip$ the inequalities
\[ 1 \leq j_u \leq i_u\ \ \ \textrm{ for all } u \in \unar \textrm{.} \]
The classes $\dj$ come from the cycles $\sigma_k^*[G_k]$ with $k=\sum_{u=1}^{r}j_u +s$.\\ \\
The cycles $\sigma_k^*[G_k]$ will arise when $k$  is such that
\[ r \leq k \leq m\ \ \textrm{, this is when } \ r \leq \sum_{u=1}^{r} j_u +s \leq \sum_{u=1}^{r} i_u + \sum_{v=1}^{s} \# A_v \textrm{.}\]
All the contributions which come from $\sigma_k^* [G_k]$ are associated to the coefficient
\[
-\frac{1}{(\sum j_u+s)!} \prod_{u=1}^{r} (j_u-1)! (-1)^{j_u-1} \ \ \ \ \ \textrm{ (a)}
\]
which appears in the sum (\ref{audessus}) .
 To choose $J_1,\ldots , J_r$, $\{ b_1 \}, \ldots ,\{ b_s \}$ amounts to choosing $r$ disjoint subsets of cardinals $j_1,j_2,\ldots,j_r$ and to order the $s$ singletons left. We count
\[ \frac{ (\sum j_u+s)!}{j_1! \ldots j_r!}\textrm{ such ways. } \ \ \textrm{ (b)} \]
To finish we have to choose an ordered partition  $\mathcal{P}$ of $\unam$ such that for all $u \in \unar$ we have $\psip [\Delta^{J_u}] = \Delta^{I_u}$, and that for each $v \in \psip \Delta^{b_v}=\Delta^{A_v}$. It amounts to choosing ordered partitions of the sets $I_u$ in $j_u$ sets. The symbol  $\lbrace^a_b \rbrace$ stands for the Stirling number of the second kind; this is the number of ways to partition a set of $a$ objects into $b$ nonempty sets\footnote{One could consult chapter 6 of \cite{gkp} to learn more about the Stirling numbers of the second kind.}. So the number of ways to partition into $b$ ordered sets is $b! \lbrace^a_b \rbrace$. In this case, we count
\[ \prod_{u=1}^{r} j_u ! \lbrace^{i_u}_{j_u}\rbrace \ \ \ \ \textrm{ (c)} \]
\noindent admissible ways to choose  $\mathcal{P}$. \\ \\
We will distinguish four cases depending on whether the sets $I_u$ are singletons (in this case ${\Phi^m}^* [\hypm]$ doesn't contribute) and whether the sets $A_v$ are singletons. \\ \\
\underline{First case: both $I_u$ and $A_v$ are singletons.} \\ \\
Such classes must come from ${\Phi^m}^*[\hypm]$, with coefficient $1$ by (\ref{margot}). So it equals  the following product when for all $u \in \unar $  we have  $i_u=1$ :
\[ \prod_{u=1}^{r} (i_u-1)! (-1)^{i_u-1} \]
\ \\ \\
\underline{Second case: the sets $I_u$ are singletons but one of the set $A_v$ is not.} \\ \\
There is as above one contribution with coefficient $1$ which comes from ${\Phi^m}^*[\hypm]$.\\
\noindent The cardinals $i_1, \ldots, i_r$ which are equal to  $1$ make the other contributions come from the cycles $\sigma_{r+s}^{*}[G_{r+s}]$. Precisely they come from classes $\dj$ with \mbox{$j_1=\ldots=j_r=1$}. We have already determined the  coefficients (a), (b) and (c) above, so the contribution we find is
\[ -\frac{(r+s)!}{(r+s)!}=-1 \textrm{.}\]
\noindent which cancels the contribution of ${\Phi^m}^*[\hypm]$. \\ \\
\underline{Third case: one of the sets $I_u$ is not a singleton but the sets $A_v$ are singletons} \\ \\
In this case, the pullback ${\Phi^m}^*[\hypm]$ does not contribute. On the other hand, we have to enumerate the contributions from the cycles $\sigma_{\sum j_u +s}^*[G_{\sum j_u +s}]$ for the integers $j_u$ which satisfy  $1\leq j_u \leq i_u$ and $\sum j_u +s < m$. The second inequality is equivalent to $\sum j_u < \sum i_u$. The sum is taken over the integers $j_u$ which verify  $1 \leq j_u \leq i_u$ and such that there exists $u \in \unar$ making $j_u$ different from $ i_u$. It is equivalent to consider the sum of all $j_u$ such that  $1 \leq j_u \leq i_u$ and to subtract the contribution corresponding to the term for which $i_u=j_u$ for each $u \in \unam$. The total contribution is then
\[
- \Big{(} \sum_{j_1=1}^{i_1} \ldots \sum_{j_r=1}^{i_r} \ \ \prod_{u=1}^{r}  (-1)^{j_u-1}  (j_u-1)!  \lbrace^{i_u}_{j_u} \rbrace    - \prod_{u=1}^{r}  (-1)^{i_u-1}  (i_u-1)!    \Big{)} \textrm{.}
\]
Notice the first sum can be factorized by
\[  \sum_{j_1=1}^{i_1} (-1)^{j_1-1} (i_1-1)! \lbrace^{i_1}_{j_1}\rbrace\]
\noindent which vanishes by theorem \ref{sommenulle}. So the class of $\di$ appears with the coefficient
\[ (-1)^{\sum i_u + r}(i_1-1)! \ldots (i_r-1)!\]
\noindent and the theorem is proven for $n=m$. \\ \\
\underline{Fourth case: one of the sets $I_u$ is not a  singleton and one of the sets  $A_v$ is not a singleton either} \\ \\
Likewise, the cycle ${\Phi^m}^*[\hypm]$ does not contribute whereas the cycles $\psip_* \dj$ from the pushdown $\sigma_{\sum j_u +s}^*[G_{\sum j_u +s}]$ contribute as soon as the  inequalities $1 \leq j_u \leq i_u$ and $\sum j_u +s < m$ hold. But we have $m=\sum i_u + \sum \# B_v$, and at least one of the cardinals $\# B_v$ is strictly greater than $1$. So the first inequality implies the second one, and the class of $\di$ appears with the coefficient
\[
 \sum_{j_1=1}^{i_1} \ldots \sum_{j_r=1}^{i_r} \ \prod_{u=1}^{r}   (-1)^{j_u-1}  (j_u-1)!  \lbrace^{i_u}_{j_u} \rbrace
\]
\noindent This coefficient vanishes as explained above. \\ \par

We used above the lemma
\begin{lem}\label{sommenulle}
For an integer $i \geq 2$ the following sum vanishes
\[ G(i)=\sum_{j=1}^{i} (-1)^{j-1}(j-1)! \lbrace^i_j \rbrace \]
\end{lem}
\begin{proof} Use the well-known identity $ \lbrace^{i}_{j} \rbrace = j \lbrace^{i-1}_{j} \rbrace + \lbrace^{i-1}_{j-1} \rbrace$ (one can consult table 250 of \cite{gkp}) and thus replace $\lbrace^{i}_{j} \rbrace$ in the definition of $G(i)$.\end{proof}

\subsection{} Now we can  simplify theorem \ref{classe_Gk} to obtain theorem \ref{th_gk_general_intro}. Notice that for integers $i_1,\ldots,i_r$ such that $\sum \ i_u \leq n$ and points  of the curve $o,o_1,\ldots,o_{n-\sum i_u}$, we have

\[ [\delta_{i_1,\ldots,i_r}+o_1+\ldots+o_{n-\sum i_u}]=[\delta_{i_1,\ldots,i_r}+(n-\sum i_u)o] \textrm{\ \ \ in } A^{g-r}(C_n) \textrm{.}
\]
 The number of ways to choose $n-\sum i_u$ distinct points in the support of a divisor $D \in \mathcal{G}$, that is, among $d$ elements, is $\binom{d}{n-\sum i_u}$. The theorem \ref{th_gk_general_intro} follows.

\section{Relations in the tautological ring modulo algebraic equivalence}\label{rel}
In this section we explain how to deduce  relations between the components $C_{(i)}$ modulo algebraic equivalence from  theorem \ref{th_gk_general_intro}.
\subsection{} By definition, all the divisors of a linear system are linearly equivalent. The morphism $u_d$ contracts the $r$-dimensional variety $\GG$ into a point of $JC$ so, by definition of the pushdown, we have
\begin{equation}\label{the_point}
 {u_d}_*[\GG]=0 \ \ \ \textrm{ in } CH^{g-r}(JC)
\end{equation}
We have already calculated the class of $\GG$ (because $\GG=G_d$)  in terms of translates of the generalized diagonals $ \delta_{i_1,\ldots,i_r}+o_1+\ldots + o_{m} =
\{ i_1x_1+\ldots+i_rx_r +o_1 + \ldots + o_{m} \ \mid \ x_i \in C \}
$. The morphism $u_d$ maps these diagonals onto translates of the variety $i_1C+ \ldots + i_rC$ and is one to one. In order to simplify the problem, we will consider relations modulo algebraic equivalence, but one could also express relations modulo the translations in the Chow ring of the Jacobian. From the definition of the Pontryagin product we have ${i_1}_*C * \ldots * {i_r}_*C  = n \ [i_1C+ \ldots + i_rC]$ where $n$ is the degree of the addition morphism
$$ \xymatrix{  i_1C \times  \ldots \times i_rC \ar[d]_n  \\ i_1C+ \ldots + i_rC    }$$
To each sequence $i_1,\ldots,i_r$ of integers we associate integers $d_1,\ldots,d_s$ such that
\[ i_1=\ldots=i_{d_1} \textrm{ , } i_{d_1} \neq i_{d_1+1} \]
\[ i_{d_1+1}=\ldots=i_{d_1+d_2} \textrm{ , } i_{d_1+d_2} \neq i_{d_1+d_2+1} \]
\[ \ldots \]
\begin{equation}\label{les_indices}
 i_{d_1+\ldots+i_{d_{s-1}+1}} = \ldots=i_{d_1+\ldots+i_{d_s}} \textrm{  and } d_1+\ldots +d_s=r
\end{equation}
and write $\gamma$ for the product
\begin{equation}\label{les_indices2}
 \gamma(d,i_1,\ldots,i_r)=d_1! \ldots d_s!
\end{equation}
With this notation, the number of antecedents of a generic point of $i_1C+ \ldots +i_rC$ is $d_1! \ldots d_s !$, so we have
\begin{equation}\label{image_ud}
 {i_1}_*C * \ldots * {i_r}_*C  = d_1! \ldots d_s! \ [i_1C+ \ldots + i_rC]\textrm{.}
\end{equation}
By (\ref{the_point}), theorem \ref{th_gk_general_intro}, and remark (\ref{image_ud}), we obtain
\begin{equation}\label{relation_algebraic1}
\sum_{1 \leq i_1 \leq \ldots \leq i_r  \leq n} \frac{(-1)^{i_1+\ldots + i_r + \epsilon}}{i_1 \ldots i_r} \binom{d}{i_1+\ldots+i_r}\frac{1}{d_1! \ldots d_s!}  \ {i_1}_*C*\ldots* {i_r}_* C \ = \ 0 \textrm{.}
\end{equation}
For all integers $i_1,\ldots,i_r$ such that $1 \leq i_1 \leq \ldots \leq i_r$ we count  $\binom{r}{d_1,\ldots,d_s}$ $r$-tuple $(j_1,\ldots,j_r)$ (by permutation of the $i_k$) so we can write
\begin{equation}\label{bob}
\sum_{i_1=1}^d \ldots \sum_{i_r=1}^{d} \frac{(-1)^{i_1+\ldots + i_r + \epsilon}}{i_1 \ldots i_r} \binom{d}{i_1+\ldots+i_r} \ {i_1}_*C*\ldots* {i_r}_* C \ = \ 0
\end{equation}
Now we use the decomposition (\ref{decomposition_C}), the action of the operator $i_*$ onto the components (this is \mbox{$i_*C_{(a)}=i^{a+2}C_{(a)}$}), and multilinearity to write
$$ {i_1}_*C*\ldots* {i_r}_* C  = \sum_{a_1=1}^{g} \ldots \sum_{a_r=1}^{g} {i_1}^{a_1+2} \ldots {i_r}^{a_r+2} C_{(a_1)}*\ldots * C_{(a_r)} \textrm{.}$$
We use the inclusions \mbox{$A^p_{(s)}*A^q_{(t)} \subset A^{p+q-g}_{(s+t)}$} and we project (\ref{bob}) onto $A^{g-r}_{(s)}(JC)$ for all $s \geq 0$ to prove theorem \ref{relations_intro}.
\subsection{}\label{subsection_study_beta} The main relation we quote is obtained in $A^{g-r}_{(s)}(JC)$. For small values of $s$, the relation we consider is trivial. For example, for $s=0$, the relation is monomial with $C_{(0)}^{*r}$. It should be trivial, otherwise we would obtain the vanishing of a small power of the theta divisor.  The first nontrivial relation is obtained for $s=d-r+1$ as stated in theorem \ref{study_beta} we prove in this section
\begin{theo}\label{study_beta}
\[ \begin{array}{lclcrc}
    \beta(d,a_1,\ldots,a_r) & =  & 0 &  \textrm{ if } &  a_1+\ldots+a_r<d-r+1 & \textrm{ and } \\ \\
     \beta(d,a_1,\ldots,a_r) &  = & (-1)^{d} \  a_1! \ldots a_r! & \textrm{ if } &   a_1+ \ldots a_r =d-r+1      &
\end{array} \]
\end{theo} \ \\
\begin{proof}We denote  $s=i_1+ \ldots +i_r$ and expand $(s-i_1 -\ldots -i_{r-1})^{a_r}$ using Newton's multinomial formula to get
\begin{equation}
\beta(d,a_1,\ldots,a_r)=\sum_{s \geq 0} (-1)^{s} \binom{d}{s} \sum_{\scriptstyle b,b_1,\ldots,b_{r-1}\atop \scriptstyle b+b_1+\ldots b_ {r-1}=a_r} (-1)^{b_1+\ldots +b_{r-1}}
\binom{a_r}{b,b_1,\ldots,b_{r-1}}
\atop s^b \ \sum_{i_1=0}^{s} \sum_{i_2=0}^{s-i_1} \sum_{i_{r-1}=0}^{s-i_1-\ldots -i_{r-1}} {i_1}^{a_1+b_1} \ldots {i_{r-1}}^{a_{r-1}+b_{r-1}}\textrm{.}
\end{equation}
\noindent By lemma \ref{sum_of_monom} below, the lower part of the sum above  is a degree $\sum_{i=1}^{r-1} a_i + \sum_{i=1}^{r-1} b_i \ + b + r-1= \sum_{i=1}^{r} a_i + r-1$ polynomial in $s$ whose leading coefficient is
\[ \frac{\prod_{i=1}^{r-1} (a_i+b_i)!}{(\sum_{i=1}^{r-1} a_i \sum_{i=1}^{r-1} b_i + r-1)! } \textrm{.}\]
The equality  $\sum_{s=0}^{d}(-1)^s (-1)^s \binom{d}{s}s^a \ = \ d! \{^a_d \} (-1)^d
$ is the  result (6.19) in  \cite{gkp}. But  the integer $\{^a_d \}$ vanishes for $a<d$ and equals $1$ when $a=d$. So when $\sum_{i=1}^{r} a_i \ < \ d-r+1$, the sum $\beta(d,a_1,\ldots,a_r)$ vanishes. In the case of equality $\sum_{i=1}^{r} a_i \ = \ d-r+1 $, it remains to compute the sum
$$ (-1)^d \ d! \sum_{\scriptstyle b,b_1,\ldots,b_{r-1}\atop \scriptstyle b+b_1+\ldots b_ {r-1}=a_r} \binom{a_r}{b,b_1,\ldots,b_{r-1}}(-1)^{b_1+\ldots+b_{r-1}} \frac{\prod_{i=1}^{r-1} (a_i+b_i)!}{(\sum_{i=1}^{r-1} a_i + \sum_{i=1}^{r-1} b_i + r-1)! } $$
\noindent which is $ \frac{\prod_{i=1}^{r}a_i!}{( d! )} $ by lemma \ref{sum_int} below. This  proves the second assertion of theorem \ref{study_beta}.
\subsection{}
\begin{lem}\label{sum_of_monom}
For all integers $0 \leq a_1$, $a_2, \ldots , a_r $ the sum
$$ \sum_{ \scriptstyle 0 \leq i_1, \ldots, i_n \atop  \scriptstyle i_1 + \ldots + i_n \leq s  } {i_1}^{a_1} \ldots {i_n}^{a_n}$$
is a polynomial of degree $a_1 + \ldots a_n +n$ in $s$. Its leading term is
\[ \frac{\prod_{i=1}^{n}a_i!}{\big{(} n +\sum_{i=1}^{n}a_i \big{)}! } \textrm{.} \]
\end{lem}
\noindent \textit{Proof: } Suppose that the lemma \ref{sum_of_monom} is proven up to rank $n$, and write the sum above with $n+1$ terms in this way:  \\
\begin{equation}\label{rackir}
\sum_{i_1=0}^{s}{i_1}^{a_1} \Big{(} \sum_{ \scriptstyle 0 \leq i_2, \ldots, i_{n+1} \atop  \scriptstyle i_2 + \ldots + i_{n+1} \leq s - i_1  } {i_2}^{a_2} \ldots {i_n}^{a_n}   \Big{)}\textrm{.}
\end{equation}
By the induction hypothesis, the term between parentheses is a degree $(a_2+\ldots+a_{n+1}+n)$ polynomial in $s-i_1$. Its leading coefficient is $\frac{a_2! \ldots a_{n+1}!}{(a_2 + \ldots + a_{n+1} + n)!}$. We will note it $P$. There exists a polynomial $Q$ whose degree is less than $a_2+\ldots+a_{n+1}+n$ and such that
$$P(s)=\frac{a_2! \ldots a_{n+1}!}{(a_2 + \ldots + a_{n+1} + n)!} s^{a_2+\ldots+a_{n+1}} +Q(s)\textrm{.}$$  \noindent So
\begin{equation}\label{jolapatate}
 \sum_{i_1=0}^{s}P(s-i_1)=\frac{a_2! \ldots a_{n+1}!}{(a_2 + \ldots + a_{n+1} + n)!} \sum_{i_1=0}^{s} (s-i_1)^{a_2+\ldots+a_{n+1}+n} {i_1}^{a_1} \atop + \sum_{i_1=0}^{s} Q(s-i_1) {i_1}^{a_1}
\end{equation}
We note $N=a_2+\ldots+a_{n+1}+n$ for brevity and develop $(s-i_1)^N$ using Newton's binomial formula. The first sum of expression (\ref{jolapatate}) becomes
\begin{equation}\label{equation2}
 \sum_{k=0}^{N} (-1)^k s^{N-k} \binom{N}{k} \Big{(} \sum_{i_1=0}^{s}{i_1}^{a_1+k} \Big{)}\textrm{.}
\end{equation}
Now use the induction hypothesis when $n=1$ to prove that the sum between parentheses in  (\ref{equation2}) is a degree $a_1+k+1$ polynomial in $s$ whose leading term is $\frac{1}{a_1+k+1}$. So (\ref{equation2})  equals
\[
\sum_{k=0}^{N}(-1)^{k} \binom{N}{k} \frac{1}{a_1+k+1}\textrm{.}
\]
\noindent The simplification of this last sum given by formula (5.51) in \cite{gkp} enables us to conclude. \end{proof}
In subsection \ref{subsection_study_beta} we have used the equality \ \\
\begin{lem}\label{sum_int}
For $r$ integers  $a_1$, $a_2,\ldots, a_r$,
\[ \sum_{\scriptstyle b,b_1,\ldots,b_{r-1}\atop \scriptstyle b+b_1+\ldots b_ {r-1}=a_r} \binom{a_r}{b,b_1,\ldots,b_{r-1}}(-1)^{b_1+\ldots+b_{r-1}} \frac{\prod_{i=1}^{r-1} (a_i+b_i)!}{(\sum_{i=1}^{r-1} a_i + \sum_{i=1}^{r-1} b_i + r-1)! } \atop = \frac{\prod_{i=1}^{r}a_i!}{(\sum_{i=1}^{r}a_i+r-1)!} \textrm{.}\]
\end{lem}
\begin{proof} Write the sum in the following way
\[
\sum_{\scriptstyle b,b_1,\ldots,b_{r-1}\atop \scriptstyle b+b_1+\ldots b_ {r-1}=a_r} \prod_{i=1}^{r-1} \Big{(} (-1)^{b_i}\frac{ (a_i+b_i)!}{b_i!}\Big{)}   \Big{(} \frac{1}{b!(\sum_{i=1}^{r} a_i + r-1 -b)!} \Big{)} \textrm{.}
\]
\noindent One could consult page 321 of \cite{gkp} for the following equalities which hold for all $i \in \{ 1,\ldots,r-1\}$:
\begin{equation}\label{bireli}
\sum_{b_i \geq 0} \frac{(-1)^{b_i}(a_i+b_i)!}{b_i!}X^{b_i} = a_i! \sum_{b_i \geq 0} \binom{a_i+b_i}{b_i}(-X)^{b_i} = \frac{a_i !}{(1+X)^{a_i+1}}\textrm{.}
\end{equation}
We also notice that
\begin{equation}
\sum_{b \geq 0} \frac{1}{b!(\sum_{i=1}^{r}a_i+(r-1)-b)!} X^{b}= \frac{1}{(\sum_{i=1}^{r} a_i + r-1)!}\sum_{b \geq 0} \binom{\sum_{i=1}^{r}a_i +r-1}{b} X^{b} \atop = \frac{(1+X)^{\sum_{i=1}^{r}a_i+r-1}}{(\sum_{i=1}^{r} a_i + r-1)!}
\end{equation}
\noindent so the sum we study is the product of $\frac{\prod_{i=1}^{r}a_i!}{(\sum_{i=1}^{r}a_i+r-1)!} $ by the coefficient of $X^{a_r}$ in the formal series
$$ \frac{1}{(1+X)^{a_1+1}}  \ldots \frac{1}{(1+X)^{{a_{r-1}}+1}} \  (1+X)^{a_1+ \ldots + a_r + r-1}\textrm{.} $$
\noindent This is the coefficient of $X^{a_r}$ in $(1+X)^{a_r}$, so it is $1$. \end{proof}
\section{Linear systems and $(2r-2)$-secant $r-2$ planes}\label{sec}
The goal of this section is to analyze the relations $C_{(i)}=0$ we can deduce from the \mbox{theorem \ref{relations_intro}}. Let us denote :
  \[
A(r,d,g) \ = \ \sum_{i=0}^{r-1} \frac{(-1)^i}{d-2r+2} \binom{i+g+r-d-2}{i} \binom{d-2r}{r-1-i} \binom{d-r+1-i}{r-i}
\]

We will prove in paragraphs \ref{vanishing_theorem_calculus} and \ref{vanishing_theorem_Ci} the :
\begin{theo}\label{vanishing_theorem}
If $C$ admits a base point free $g^r_d$ and if $A(r,d,g) \neq 0$, then we have  \\
\[   C_{(i)}=0 \textrm{ in } A^{g-1}_{(i)} \textrm{ for } i \geq d-2r+1.\]
\end{theo}
But we will explain  in \ref{vanishing_theorem_cvg} how it can be deduced from the theorem of Colombo and van Geemen because $A(r,d,g)$ is the number of projective spaces of dimension $r-2$ which cut the curve in $2r-2$ points.\footnote{It could suggest the existence of a partial converse to the theorem of Colombo and van Geemen.}\\ \\

\subsection{}\label{vanishing_theorem_calculus}
By  theorem \ref{study_beta}, the first nontrivial  relation is obtained in $A^{g-d+2r-1}_{(d-2r+1)}(JC)$:
\begin{equation}\label{gaspard}
\sum_{\scriptstyle 1 \leq a_1 \leq \ldots \leq a_r \atop \scriptstyle a_1 +\ldots a_r =d-2r+1} \frac{(a_1+1)! \ldots (a_r+1)!}{\gamma(a_1,\ldots,a_r)} \ C_{(a_1)} \ldots C_{(a_r)} \textrm{.}
\end{equation}
Now multiply (\ref{gaspard}) by $C_{(0)}^{*(g+r-d-2)}$ and order the terms with respect to the powers of $C_{(0)}$. We introduce $t$ which is the number of integers $a_i$ which are equals to $1$. So we have \mbox{$a_1=\ldots=a_t=1$} and $a_{t+1}>1$. By  definition of $\gamma$ we have $\gamma(a_1,\ldots,a_r)=t! \ \gamma(a_{t+1},\ldots,a_r)$, and
\begin{equation}\label{duke}
\sum_{t=0}^{r-1} \sum_{ \scriptstyle 1 < a_{t+1} \leq \ldots \leq a_r \atop \scriptstyle a_{t+1} + \ldots +a_r=d-2r+1-t} \frac{(a_{t+1}+1)!\ldots (a_r+1)!}{t! \gamma(a_{t+1},\ldots,a_r)} C_{(0)}^{*(t+g+r-2-d)}*C_{(a_{t+1})} \ldots C_{(a_r)} \ = \ 0 \textrm{.}
\end{equation}
By Fourier transform, $R^{g-1}_{(d-2r+1)}$ is isomorphic to $R^{d-2r+2}_{(d-2r+1)}$. This  subspace  is spanned by $C_{(d-2r+1)}$, so its dimension is bounded by $1$. The relations Polishchuk obtained  enable us to express $C_{(0)}^{*(t+g+r-2-d)}C_{(a_{t+1})} * \ldots * C_{(a_r)}$ in terms of $C_{(0)}^{*(g+2r-d-3)}*C_{(d-2r+1)}$. The corollary of lemma 0.3 in \cite{poli} gives
\[
C_{(0)}^{*(g-1-\sum_{i=1}^{k} n_i)}C_{(n_1)} \ldots C_{(n_k)} \ = \ (-1)^{k-1}\frac{(g-1-\sum_{i=1}^{k} n_i))!}{(g+k-2-\sum n_i)!}\frac{(\sum n_i)!}{n_1! \ldots n_k!} C_{(0)}^{*(g+k-2-\sum n_i)} C_{(1-k+\sum n_i)}
\]
Thus (\ref{duke}) becomes a monomial relation whose coefficient is
\[
\sum_{t=0}^{r-1} (-1)^{r+t} \frac{(t+g+r-d-2)!(d-r+1-t)!}{t! (r-t)!}  \Big{(}\sum_{ \scriptstyle 1 < a_{t+1} \leq \ldots \leq a_r \atop \scriptstyle a_{t+1} + \ldots +a_r=d-2r+1-t} \frac{(r-t)!}{\gamma(a_{t+1},\ldots,a_r)}\Big{)}\textrm{.}
\]
We note that $\gamma(a_{t+1},\ldots,a_r)=\gamma(a_{t+1}-1,\ldots,a_r-1)$ and reindex the sum to obtain
\[
\sum_{ \scriptstyle 1 \leq a_{t+1} \leq \ldots \leq a_r \atop \scriptstyle a_{t+1} + \ldots +a_r=d-2r+1-t} \frac{(r-t)!}{\gamma(a_{t+1},\ldots,a_r)} = \sum_{ \scriptstyle 0 \leq b_{1} \leq \ldots \leq b_{r-t} \atop \scriptstyle b_1 + \ldots +b_{r-t}=d-3r+1+t} \frac{(r-t)!}{\gamma(b_1,\ldots,b_{r-t})} \textrm{.}
\]
As $\frac{(r-t)!}{\gamma(b_1,\ldots,b_{r-t})}$ is the number of $(r-t)$-tuples we can associate to the integers $\{ b_1\} $, $\{ b_2\},\ldots, \{ b_{r-t}\} $, the sum corresponds to the number of positive integers whose sum is  $d-3r+1+t$. If we consider the formal series $\frac{1}{1-X}=1+X+X^2+\ldots$, the integer we are looking  for appears as the coefficient of $X^{d-3r+1+t}$ in $\Big{(} \frac{1}{1-X} \Big{)}^{r-t}$, that is
\[
\Big{(}\sum_{t=0}^{r-1} (-1)^{r+t} \frac{(t+g+r-d-2)!(d-r+1-t)!}{(r-t)!t!} \binom{d-2r}{d-3r+1+t} \Big{)} \ \ C_{(0)}^{*(g+2r-d-3)} * C_{(d-2r+1)} \ =  \ 0 \textrm{.}
\]
We factor by $(g+r-d-2)!(d-2r+2)!$ to  make the coefficient $A(r,d,g)$ below appear. When it is not zero, we have $C_{(d-2r+1)} \ = \ 0$ modulo algebraic equivalence. Theorem \ref{vanish_theo} finishes the proof.
\subsection{}\label{vanishing_theorem_Ci}
The following remark is an easy consequence of the relations obtained by Alexander Polishchuk. Giambattista Marini has also noticed and proven it in \cite{marini}.
\begin{theo}\label{vanish_theo}
Let $i$ be a positive integer. If $C_{(i)}=0$ in $A^{g-1}_{(i)}$ then $C_{(i+1)}=0$ in $A^{g-1}_{(i+1)}$.
\end{theo}
\begin{proof}If $g \in \{ 2,3\}$, the only components which arise are $C_{(0)}$ and $C_{(1)}$, and the cycle $C_{(0)}$ never vanishes because $-\fo(C_{(0)})$ is the class of the theta divisor (see \cite{arnaudc}). Now suppose that $g \geq 4$. As the gonality of the curve is less than $g$, apply the theorem of Colombo and van Geemen to get $C_{(g-2)}=0$ and the implication $C_{(g-3)}=0 \ \ \Rightarrow \ \ C_{(g-2)}=0$. Now suppose that $i \leq g-4$. The point is that $R^{g-1}_{(i+1)}$, which is spanned by $C_{(i+1)}$, is isomorphic to $R^{i+1}_{(i)}$, which is spanned by $C_{(0)}^{*(g-i-2)}*C_{(i+1)}$ by  corollary 0.3 in \cite{poli}. Now consider the cycle $C_{(0)}*C_{(1)}*C_{(i)}$ which is proportional to $C_{(0)}^{*(g-i-2)}*C_{(i+1)}$ by the same corollary.\end{proof}
\subsection{}\label{vanishing_theorem_cvg}
The Grassmanian  $\mathbb{G}(r-1,r+1)$ which parametrises the $(r-2)$-dimensional linear  subspaces of $\Pr$ is a $(2r-2)$-fold. The subspaces containing a point of the curve correspond to a codimension $1$ subvariety of $\mathbb{G}(r-1,r+1)$. So one can expect to find a finite number of  $\PV^{r-2}$ which contains $2r-2$ points of the curve. Actually,  the formula of Castelnuovo computes this number when it is finite:
\begin{equation}\label{def_crdg}
B(r,d,g)=\sum_{i=0}^{r-1}\frac{(-1)^i}{r-i}\binom{d-r-i+1}{r-1-i} \binom{d-r-i}{r-1-i} \binom{g}{i}\textrm{.}
\end{equation}
One could consult \cite{ELMS} for a proof. When $B(r,d,g)$ is nonzero, consider $2r-2$ such points and then the hyperplanes which contain these points. The linear system we obtain in this way is a $g^1_{d-2r+2}$. Applying the theorem of Colombo and van Geemen we obtain an analogue of theorem \ref{vanishing_theorem} with the condition on $B(r,d,g)$ : when $B(r,d,g) \neq 0$ then $C_{(i)}=0$ in $A^{g-1}_{(i)}$ for $i \geq d-2r+1$. But the theorems are the same as
\begin{theo}
For each integer $ r \geq 1$ we have the equality in $\Q[d,g]$:
\[ A(r,d,g)=B(r,d,g)\]
\end{theo}

We use the algorithm of Zeilberger to prove that $A(r,d,g)=B(r,d,g)$. We apply it thanks to the software Mathematica and the package ZW that Peter Paule and Markus Schorne have developed. It enables us to compute a linear recurrence relation
\[ Q_2(r,d,g)X(r+2,d,g) + Q_1(r,d,g)X(r+1,d,g) + Q_0(r,d,g)X(r,d,g)  \ = \ 0 \]
with coefficients $Q_i \in \Q[r,d,g]$ that both $A$ and $B$ satisfy. We conclude because the first terms  are the same: $A(1,d,g)=B(1,d,g)$ and $A(2,d,g)=B(2,d,g)$. To learn more about Zeilberger's algorithm,  one can consult \cite{a_equals_b}. One can also find an explanation of the ZW package in \cite{zw}.

\section{Relations in the tautological ring for pencils, nets and webs. }\label{app}
Theorem \ref{relations_intro} gives immediately the theorem of Colombo and van Geemen (see paragraph \ref{cvg_theorem}). We have seen in the last section that the relations $C_{(i)}=0$ we can deduce from the theorem \ref{relations_intro} could also be deduced  from their theorem. Do we obtain new relations modulo algebraic equivalence? Yes, we give for plane and space curves of low genus the new relations in table (\ref{tab_plane_curves}) and (\ref{tab_space_curves}).

\subsection{}\label{cvg_theorem}
If $C$ admits a base point free $g^1_d$ then theorem \ref{relations_intro} gives, for $0 \leq a  \leq g-1$:
\[ \beta(d,a+1) \ C_{(a)} \ = \ 0 \ \ \ \textrm{ in } A^{g-1}_{(a)}\textrm{.}    \]
The coefficients $\beta(d,a)$ express in terms of the Stirling numbers of the second kind. One can consult formula  (6.19) of \cite{gkp}:
\[
\beta(d,a)=\sum_{i=1}^{d}(-1)^i (-1)^i \binom{d}{i}i^a=d! \{^a_d \} (-1)^d\textrm{.}
\]
Lastly, $\{^{a+1}_{ \ d } \}$ is nonzero if and only if $0 \leq d  \leq a+1$, so we have  $C_{(a)}=0$ for $a \geq d-1$. This is Colombo and van Geemen's result.
\subsection{}
The integer $A(2,d,g)=\frac{(d-1)(d-2)}{2}-g$ is the well-known number of nodes  for a plane curve of degree $d$ and genus $g$. So for every plane curve the integer $A(2,d,g)$ is positive. Conversely, when the integers $d$ and $g$ satisfy $g \geq 0$ and $\frac{\sqrt{8g+1}+3}{2} \leq d$, there exists a plane curve with such invariants. As explained in the last section, when the curve is singular one can consider the line passing through a singular point to obtain a $g^1_{d-2}$. Let us denote $d'=d-2$ in this case. When the curve is smooth, it only admits a $g^1_{d-1}$. We denote $d'=d-1$ in this case. We indicate in the following table the relations\footnote{Keep in mind that Polishchuk used the cycles $p_i=\fo(C_{(i-1)})$ to describe $R$. So to translate the relations, one has to change $C_{(i)}$ into $p_i$ and the Pontryagin products into intersection products.} deduced from  theorem \ref{relations_intro}. The new relations (those which cannot be deduced from $I_g$ and from Colombo and van Geemen's theorem  applied to the $g^1_{d'}$ the curve admits)\footnote{More precisely, to define relations one should consider the kernel  $\mathcal{R}$ of the map \mbox {$\pi : \C[C_{(0)},\ldots,C_{(0)}] \longrightarrow R$} which maps the indeterminate $C_{(i)}$ to the cycle $C_{(i)}$. One could also define (with Polishchuk's formula) a Fourier transform on $\C[C_{(0)},\ldots,C_{(0)}]$ compatible with $\fo$. In this terms, a new relation is an element of $\mathcal{R}$ which is not in the smallest ideal containing $I_g$, $C_{(d'-1)}$ and stable under Fourier transform.} are written \textbf{in bold type}. \\ \par
\begin{table}
\centering
\caption{Relations deduced from Theorem \ref{relations_intro} for plane curves of genus $g \leq 9$}
\begin{tabular}{| l || c || c || r || r |}\hline
Genus & $g^2_d$  & $g^1_{d'}$ deduced & Relations deduced from the $g^2_d$ & Consequences \\
              &                     &  from the $g^2_d$       &   with Th. \ref{relations_intro}                                                                 &   \\ \hline \hline
$g=5$ & $g^2_5$ & $g^1_3$ & $3C_{(0)} *C_{(2)} +C_{(1)}^{*2}=0$&   \\ \hline
$g=6$ & $g^2_5$& $g^1_4$ & $\boldsymbol{3C_{(0)}  C_{(2)}+C_{(1)}^2=0}$ & \\ \hline
$g=7$ &$g^2_6$ &$g^1_4$ & $2C_{(0)}  *  C_{(3)}+C_{(1)}  *  C_{(2)}=0$ & \\ \hline
$g=8$  & $g^2_7$&  $g^1_5$ & $\boldsymbol{8C_{(1)}*C_{(3)}+3C_{(2)}^{*2}=0}$ &
 $\boldsymbol{C_{(1)}*C_{(3)}=C_{(2)}^{*2}=0}$                                                \\ \cline{2-5}

	       & $g^2_6$  & $g^1_4$            &   $\boldsymbol{C_{(1)}*C_{(2)}=0}$   &                                                    $\boldsymbol{C_{(0)}*C_{(1)}*C_{(2)}=0}$                                               \\
	       & & & & $\boldsymbol{C_{(1)}^{*3}=0}$ \\ \hline
& &  & &$\boldsymbol{C_{(1)}*C_{(3)}=C_{(2)}^{*2}=0}$  \\
& $g^2_7$& $g^1_5$  & $\boldsymbol{8C_{(1)}*C_{(3)}+3C_{(2)}^{*2}=0}$ & $\boldsymbol{C_{(0)}*C_{(1)}*C_{(3)}=C_{(0)}*C_{(2)}^{*2}=0}$ \\
$g=9$& & & & $\boldsymbol{C_{(1)}^{*2}*C_{(2)}=0}$\\ \cline{2-5}
 & & &  & $\boldsymbol{C_{(0)}*C_{(1)}*C_{(2)}=0}$ \\
 &$g^2_6$ & $g^1_4$& $\boldsymbol{C_{(1)}*C_{(2)}=0}$ &$\boldsymbol{C_{(0)}^{*2}*C_{(1)}*C_{(2)}=0}$ \\
 & & & &$\boldsymbol{C_{(0)}*C_{(1)}^{*3}=0}$ \\ \hline
 \end{tabular}
\label{tab_plane_curves}
\end{table}

\noindent \textit{Remark: }For $g$ sufficiently large, one can always find $d$ such that it exists a genus $g$ curve with a $g^2_d$ giving us new relations. Let us sketch the proof. The first non-trivial relation given by the theorem \ref{relations_intro} is
\begin{equation}\label{new_relation}
 \sum_{\scriptstyle 1\leq a,b \atop \scriptstyle a+b=d-3} \ (a+1)!(b+1)! \ C_{(a)}*C_{(b)} \ = \ 0 \ \ \ \textrm{ in } R^{g-2}_{(d-3)}\textrm{.}
\end{equation}
\noindent  The relations of $I_g$ correspond to subspaces $R^{d+i}_{(i)}$ for $0\leq d \leq \frac{g}{2}-1$ and $d+1 \leq i \leq g-d-1$ (one has to see how the relations of the theorem 0.1 of \cite{poli} are parametrized). So when the inequality $d<\frac{g}{2}+2$ holds, there is no relation of $I_g$ which belongs to $R^{g-2}_{(d-3)}$. Lately, if $d'=d-1$, there is no relation deduced from the existence of the $g^1_{d'}$ in $R^{g-2}_{(d-3)}$. If $d'=d-2$, the only relation we can deduce from the $g^1_{d'}$ in $R^{g-2}_{(d-3)}$ is $C_{(0)}*C_{(d-3)}=0$ which is different from (\ref{new_relation}).
\subsection{}
The integer $A(3,d,g)$ is the number of quadrisecants to the space curve $\Phi(C)$. Let us recall Cayley's formula
\begin{equation}\label{quadrisecantes}
A(3,d,g)= \frac{(d-2)(d-3)^2(d-4)}{12} \ - \ \frac{g(d^2-7d+13-g)}{2}\textrm{.}
\end{equation}
When $A(3,d,g)$ is non-zero, the curve admits a $g^1_{d-4}$ and so $C_{(d-5)}$ is zero modulo algebraic equivalence. We note $d'=d-4$ in this case. Otherwise, we deduce from the existence of a $g^1_{d-3}$ the vanishing of $C_{(d-4)}$. We note  $d'=d-3$. We give in the following table the results we obtain by theorem \ref{relations_intro} for space curves of low genus. As above, a new relation is written \textbf{in bold type}. It means the relation cannot be deduced from the theorem of Colombo and van Geemen applied to the $g^1_{d'}$ and from the relations of $I_g$.
\begin{table}
\center
\caption{Relations deduced from Theorem \ref{relations_intro} for space curves of genus $g \leq 9$}
\begin{tabular}{| l || c || c || r || r |}\hline
Genus & $g^3_d$  & $g^1_{d'}$ deduced & Relations deduced from $g^3_d$ & Consequences \\
              &                     &  from $g^3_d$       &   with Th. \ref{relations_intro}                                                                  &   \\ \hline \hline
$g=6$ & $g^3_7$ & $g^1_3$ & $C_{(0)}*C_{(1)}^{*2}=0$ & \\ \hline
$g=7$&$g^3_8$ &$g^1_4$ &$9C_{(0)}*C_{(1)}*C_{(2)}+C_{(1)}^{*3}=0$ & \\ \hline
$g=8$&$g^3_8$ &$g^1_4$ &$\boldsymbol{9C_{(0)}*C_{(1)}*C_{(2)}+C_{(1)}^{*3}=0}$ &$\boldsymbol{C_{(1)}*C_{(2)}=0}$ \\
& & & & $\boldsymbol{C_{(1)}^{*3}=C_{(0)}*C_{(1)}*C_{(2)}=0}$ \\ \hline
& & & $\boldsymbol{8C_{(0)}*C_{(1)}*C_{(3)}}$ & $\boldsymbol{C_{(1)}*C_{(3)}=C_{(2)}^{*2}=0} $  \\
& $g^3_9$& $g^1_5$ &$\boldsymbol{+3C_{(0)}*C_{(2)}^{*2}}$ & $\boldsymbol{C_{(1)}^{*2}*C_{(2)}=0}$\\
& & & $\boldsymbol{+2C_{(1)}^{*2}*C_{(2)}=0}$ & $\boldsymbol{C_{(0)}*C_{(2)}^{*2}=0}$\\
& & & & $\boldsymbol{C_{(0)}*C_{(1)}*C_{(3)}=0}$ \\ \cline{2-5}
$g=9$ & &  &  & $ \boldsymbol{C_{(0)}^{*m}*C_{(1)}*C_{(2)}=0}$   \\
& $g^3_8$  & $g^1_4$& $\boldsymbol{9C_{(0)}*C_{(1)}*C_{(2)}+C_{(1)}^{*3}=0}$ & $\boldsymbol{\textbf{ for } m \in \{ 0 \ldots 2\}}$ \\
& & & & $\boldsymbol{C_{(0)}^{*n}*C_{(1)}^{*3}=0}$ \\
& & & & $\boldsymbol{ \textbf{ for } n \in \{ 0 ,1\}}$ \\ \hline
\end{tabular}
\label{tab_space_curves}
\end{table}
\noindent As above, for $g$ sufficiently large there exists a genus $g$ space curve for which theorem \ref{relations_intro} gives us new relations.
\bibliographystyle{amsalpha}
\bibliography{biblio_2006.bib}
\end{document}